\documentclass{article}
\usepackage{authblk}

\usepackage[backend=bibtex,doi=false,eprint=false,firstinits=true,isbn=false,style=numeric-comp,url=false]{biblatex}
\bibliography{conti.bib}
\DeclareRedundantLanguages{english,german,french}{english,german,ngerman,french}

\usepackage[ngerman,english]{babel}
\usepackage[latin1]{inputenc}
\usepackage{amsfonts,amsmath,amsthm}
\usepackage{mathabx}
\usepackage{cases}
\usepackage{xfrac}
\usepackage{enumerate}
\usepackage{fancyhdr}
\usepackage{color}
\usepackage[colorlinks]{hyperref}
\definecolor{blue75}{rgb}{0,0,.75}
\definecolor{green75}{rgb}{0,.75,0}
\hypersetup{colorlinks=true, urlcolor=blue75,linkcolor=blue75,citecolor=green75,pdfstartview=FitB,bookmarksopen=true,bookmarksopenlevel=1}
\usepackage{csquotes}
\usepackage[a4paper, left=2.5cm, right=2.5cm, top=2.5cm,bottom=2cm]{geometry}
\usepackage{constants}
\newcommand{\parenthezises}[1]{\arabic{#1}}
\newconstantfamily{C}{
symbol=C,
format=\parenthezises,
reset={section}
}
\newconstantfamily{M}{
symbol=M,
format=\parenthezises,
reset={section}
}
\allowdisplaybreaks
\begin{document}
\newcommand{\R}{\mathbb{R}}
\newcommand{\N}{\mathbb{N}}
\def\dist{\operatorname{dist}}
\def\ess{\operatorname{ess}}
\def\osc{\operatorname{osc}}
\def\sign{\operatorname{sign}}
\def\supp{\operatorname{supp}}
\newcommand{\LOne}{L^{1}(\Omega)}
\newcommand{\LTwo}{L^{2}(\Omega)}
\newcommand{\Lq}{L^{q}(\Omega)}
\newcommand{\Lp}{L^{p}(\Omega)}
\newcommand{\LInf}{L^{\infty}(\Omega)}
\newcommand{\HOneO}{H^{1}_{0}(\Omega)}
\newcommand{\HTwoO}{H^{2}_{0}(\Omega)}
\newcommand{\HOne}{H^{1}(\Omega)}
\newcommand{\HTwo}{H^{2}(\Omega)}
\newcommand{\HmOne}{H^{-1}(\Omega)}
\newcommand{\HmTwo}{H^{-2}(\Omega)}

\newcommand{\qn}{H^2_{loc}\left(\R^+_0,H^{-\frac{1}{2}}(\Gamma)\right)}
\newtheorem{Theorem}{Theorem}[section]
\newtheorem{Assumption}[Theorem]{Assumptions}
\newtheorem{Corollary}[Theorem]{Corollary}
\newtheorem{Convention}[Theorem]{Convention}
\newtheorem{Definition}[Theorem]{Definition}
\newtheorem{Example}[Theorem]{Example}
\newtheorem{Lemma}[Theorem]{Lemma}
\newtheorem{Notation}[Theorem]{Notation}
\newtheorem{Remark}[Theorem]{Remark}
\numberwithin{equation}{section}
\pagestyle{plain}
\renewcommand\Affilfont{\itshape\small}
\title{On a new transformation for generalised porous medium equations: from weak solutions to classical}
\author{Anna Zhigun}
\affil{Technische Universität Kaiserslautern, Felix-Klein-Zentrum für Mathematik\\ Paul-Ehrlich-Str. 31, 67663 Kaiserslautern, Germany\\
  e-mail: {zhigun@mathematik.uni-kl.de}}
\date{}
\maketitle
\begin{abstract}
It is well-known that solutions for generalised porous medium equations are, in general, only H\"older continuous. In this note, we propose a new variable substitution for such equations which transforms weak solutions into classical.
\end{abstract}

\section{Introduction}
Let us consider a generalised porous medium equation in the form
\begin{subequations}\label{PMD}
\begin{alignat}{2}
 &\partial_t u=a(u)\Delta u+f(u)&&\text{ in }(0,T)\times\Omega=:Q_T,\label{PM1}\\
 &u=u_{\Gamma}&&\text{ in }(0,T)\times\Gamma,\\
 &u=u_{0}&&\text{ in }\Omega
\end{alignat}
\end{subequations}
for some $T>0$ and a bounded domain $\Omega\subset\R^d$, $d\in\N$, with a smooth boundary $\Gamma$. We assume the diffusion coefficient $a$ to be strictly positive but for $a(0)=0$. Below, we shell give detailed assumptions on $a$, the reaction term $f$, the boundary  data $u_{\Gamma}$ and the initial data $u_0$ (see {\it Assumptions~\ref{A1}}).
If, for example, $a(u)=m |u|^{1-\sfrac{1}{m}}$ for some $m>1$, then equation \eqref{PM1} can be obtained from the standard  porous medium equation 
\begin{align}
 &\partial_t s=\Delta (s|s|^{m-1})+g(s)\text{ in }Q_T,\label{PM}
\end{align}
where
\begin{align}
  s:=u|u|^{\sfrac{1}{m}-1},\ g(s):=\sfrac{1}{m}|s|^{1-m}f(s|s|^{m-1}).\nonumber
\end{align}
It is well-known that  solutions of  \eqref{PM} are, in general, only weak solutions if $s\not\equiv0$ and is not strictly separated from zero. In particular, $\partial_t s$ and $\Delta (s|s|^{m-1})$ are generally not even uniformly bounded.
Still, it is well-understood \cite{AroBen1979,DiBen83,CaffFri1980,Ziemer1982} that, under reasonable assumptions on $g$, the solution $s$ is at least H\"older continuous in $\overline{Q}_T$. The same is all the more true for $u=s|s|^{m-1}$. 
In this note, we show how the information on the H\"older continuity of solution $u$ can be used in order to transform equation \eqref{PM1} into a generalised porous medium equation  with a {\it classical} solution $v$. Our construction is based on a new variable substitution $v=V(u)$ by means of a smooth and strictly increasing function $V$ which depends only upon $T$, the structure of $\Gamma$  and some H\"older exponents and norms of $a$, $f$ and $u$. 
Thus, although $u$ itself is only H\"older continuous, it can be reconstructed from a regular
solution $v$.

This note is organised as follows: in  {\it Section~\eqref{Pr}} we state our assumptions and result, which we then prove in  {\it Section~\eqref{Proof}}. Finally, in {\it Section~\eqref{Example}}, we give an explicit transformation for a solution of the classical homogeneous porous medium equation.
\section{Problem setting and main result}\label{Pr}
In this note, we are interested in a transformation of \eqref{PMB} which leads to a classical solution. The latter means that the resulting equation should hold in some H\"older space.
We thus assume the reader to be familiar with the standard and anisotropic H\"older spaces and their standard properties. 
We make the following structural assumptions upon the parameters and a solution of problem \eqref{PMB}:
\begin{Assumption}\label{A1}~
\begin{enumerate}
 \item The diffusion coefficient $a:\R\rightarrow\R_0^+$ has the properties
 \begin{enumerate}
  \item $a(0)=0$, $a>0$ in $\R\backslash\{0\}$;
  \item $a$ is increasing on $\R^+_0$ and decreasing on $\R^-_0$;
  \item $a\in C^{(\alpha_a)}(\R\backslash\{0\})$ for some exponent $\alpha_a\in(0,1)$;
 \end{enumerate}
\item The reaction term $f:\R\rightarrow\R$ has the property: $f\in C^{(\alpha_f)}(\R\backslash\{0\})$ for some exponent $\alpha_f\in(0,1)$;
\item The function $u:\overline{Q_T}\rightarrow\R$ has the properties
\begin{enumerate}
 \item $u$ is a solution to initial-boundary value problem \eqref{PMD} in a weak \cite{Vazquez} sense;
 \item $u\in C^{\left(\sfrac{\alpha_u}{2},\alpha_u\right)}(\overline{Q_T})$ for some exponent $\alpha_u\in(0,1)$;
 \item $u_{\Gamma}\in C^{\left(1+\sfrac{\alpha}{2},2+\alpha\right)}(\{u_{\Gamma}\neq0\})$, $u_0\in C^{(2+\alpha)}(\{u_0\neq0\})$ for $\alpha:=\min\{\alpha_a,\alpha_f\}\alpha_u$, and the  compatibility condition of the first order \cite[Chapter IV, \textsection 5]{LSU} holds in the following sense: 
 \begin{align}
  u_{\Gamma}(0,\cdot)=u_0,\ a(u_0)\Delta u_0+f(u_0)=\partial_t u_{\Gamma}(0,\cdot)\text{ in }\Gamma\cap\{u_0\neq0\}.\nonumber
 \end{align}
\end{enumerate}
\end{enumerate}
\end{Assumption}
Under {\it Assumptions~\ref{A1}}, we prove the following theorem:
\begin{Theorem}\label{maint}
 Let $\Omega\subset\R^d$, $d\in\N$, be a bounded domain with a smooth boundary $\Gamma$ and let {\it Assumptions~\ref{A1}} hold. There exist: 
 \begin{enumerate}
  \item a strictly increasing function $V\in W^{1,\infty}\left[-||u||_{C(\overline{Q_T})},||u||_{C(\overline{Q_T})}\right]$ with $V(0)=0$,
  \item a function $\bar{f}_u\in C^{\left(\sfrac{\alpha}{2},\alpha\right)}(\overline{Q_T})$,
 \end{enumerate}
such that $v:=V(u)$ is a classical solution of a generalised porous medium equation:
 \begin{align}
  &\partial_t v=\nabla\cdot(a\circ U(v)\nabla v)+\bar{f}_u\text{ in } C^{\left(\sfrac{\alpha}{2},\alpha\right)}(\overline{Q_T}).\label{eq3}
 \end{align}
Here, $U$ denotes the inverse of $V$.
Function $V$ depends only upon $T$, the structure of $\Gamma$, the H\"older exponents $\alpha_a,\alpha_f$ and $\alpha_u$ and the function
\begin{align}
 &\psi:\left[-||u||_{C(\overline{Q_T})},||u||_{C(\overline{Q_T})}\right]\backslash\{0\}\rightarrow\R^+_0,\text{ and for all }k\in\left(0,||u||_{C(\overline{Q_T})}\right]\nonumber\\ &\psi(k):=\max\left\{||a||_{C^{(\alpha_a)}\left[k,||u||_{C(\overline{Q_T})}\right]},||f||_{C^{(\alpha_f)}\left[k,||u||_{C(\overline{Q_T})}\right]},||u||_{C^{\left(\sfrac{\alpha_u}{2},\alpha_u\right)}(\{u\geq k\})},\right.\\
 &\left.\quad\quad\quad\quad\quad\quad ||u_{\Gamma}||_{C^{\left(1+\sfrac{\alpha}{2},2+\alpha\right)}(\{u_{\Gamma}\geq k\})},||u_0||_{C^{(2+\alpha)}(\{u_0\geq k\})}\right\},\nonumber\\
 &\psi(-k):=\max\left\{||a||_{C^{(\alpha_a)}\left[-||u||_{C(\overline{Q_T})},-k\right]},||f||_{C^{(\alpha_f)}\left[-||u||_{C(\overline{Q_T})},-k\right]},||u||_{C^{\left(\sfrac{\alpha_u}{2},\alpha_u\right)}(\{u\leq-k\})},\right.\\
 &\left.\quad\quad\quad\quad\quad\quad\ \ ||u_{\Gamma}||_{C^{\left(1+\sfrac{\alpha}{2},2+\alpha\right)}(\{u_{\Gamma}\leq-k\})},||u_0||_{C^{(2+\alpha)}(\{u_0\leq-k\})}\right\}.\nonumber
\end{align}
\end{Theorem}
\section{Proof of {\it Theorem~\ref{maint}}}\label{Proof}
\begin{Notation}
 To shorten the notation, we make the following convention: if a quantity depends only upon such parameters as $T$, the structure of $\Gamma$, the H\"older exponents $\alpha_a,\alpha_f$ and $\alpha_u$ and the function $\psi$, we say that it depends only upon the parameters of the problem.
\end{Notation}
We need some preliminary work in order to define the desired transformation $V$.
Let the {\it Assumptions~\ref{A1}} hold. Then, it follows that
\begin{align}
 a(u)\in C^{\left(\sfrac{\alpha_a\alpha_u}{2},\alpha_a\alpha_u\right)}(\{u\neq0\}),\ f(u)\in C^{\left(\sfrac{\alpha_f\alpha_u}{2},\alpha_f\alpha_u\right)}(\{u\neq0\}).\nonumber
\end{align}
Moreover, since $u$ is continuous on $\{u\neq0\}$, the sets $\{u\geq k\}$ ($\{u\leq-k\}$) and $\{u>k\}$ ($\{u<-k\}$) are for all $k\in\left(0,||u||_{C(\overline{Q_T})}\right)$ compact and relatively (with respect to $\overline{Q_T}$) open, respectively. Thus, each set $\{u\geq k\}$ ($\{u\leq-k\}$) can be covered by a finite number of relatively (with respect to $\overline{Q_T}$) open cylinders
contained in $\{u>\sfrac{k}{2}\}$ ($\{u<-\sfrac{k}{2}\}$). In tern, for each such cylinder $Q$, there exists a number $\delta_k>0$ such that the cylinder $Q_{\delta_k}:=\{x\in Q_T|\ \dist(x,Q)<\delta_k\}$, where $\dist(x,Q):=\max\{|x-y|\ | y\in Q_T\}$, is contained in $\{u>\sfrac{k}{2}\}$ ($\{u<-\sfrac{k}{2}\}$) as well. Equation \eqref{PM1} is non-degenerate on $\{u>\sfrac{k}{2}\}$ ($\{u<-\sfrac{k}{2}\}$). Moreover, the boundary trace $u_{\Gamma}$ and the initial value $u_0$ are regular and compatible on these sets due to {\it Assumptions~\ref{A1} 3(d)}.  Therefore, we can apply a standard result on the local regularity of linear parabolic equations with H\"older continuous coefficients, Theorem 10.1 from \cite[Chapter IV, \textsection 10]{LSU}, to the cylinder $Q$ as subcylinder of $Q_{\delta_k}$. Since $\{u\geq k\}$ ($\{u\leq-k\}$) is covered by a finite number of cylinders of this type, the result of that theorem can be interpreted in the following way: there exists a function
\begin{align}
 \varphi:\left[-||u||_{C(\overline{Q_T})},||u||_{C(\overline{Q_T})}\right]\rightarrow\R^+_0\nonumber
\end{align}
with the properties
\begin{enumerate}
 \item  $\varphi(0)=0$, $\varphi>0$ in $\left[-||u||_{C(\overline{Q_T})},||u||_{C(\overline{Q_T})}\right]\backslash\{0\}$;
 \item $\varphi$ is increasing on $\left[-||u||_{C(\overline{Q_T})},0\right]$ and decreasing on $\left[0,||u||_{C(\overline{Q_T})}\right]$;
 \item for $\alpha=\min\{\alpha_a,\alpha_f\}\alpha_u$, it holds for all $k\in\left(0,||u||_{C(\overline{Q_T})}\right]$ that 
 \begin{align}
 &||u||_{C^{\left(1+\sfrac{\alpha}{2},2+\alpha\right)}(\{u\geq k\})}\leq \varphi^{-1}(k),\nonumber\\
 &||u||_{C^{\left(1+\sfrac{\alpha}{2},2+\alpha\right)}(\{u\leq -k\})}\leq \varphi^{-1}(-k).\nonumber
\end{align}
 \item $\varphi$ depends only upon the parameters of the problem.
\end{enumerate}
Without loss of generality, we may also assume  that
\begin{align}
 \varphi(k)\leq a^2(k)\text{ for all }k\in\left[-||u||_{C(\overline{Q_T})},||u||_{C(\overline{Q_T})}\right].\label{sm}
\end{align}
Otherwise, we replace $\varphi$ by $\min\{\varphi,a^2\}$.
Using $\varphi$, we construct yet another function
\begin{align}
&\Phi:\left[-||u||_{C(\overline{Q_T})},||u||_{C(\overline{Q_T})}\right]\rightarrow\R,\ \Phi(k):=\int_0^k\int_0^x\left(\int_0^y\varphi\left(\sfrac{z}{2}\right)\,dz\right)^2\,dy\,dx.\nonumber
\end{align}
It is obvious that $\Phi\in W^{3,\infty}\left[-||u||_{C(\overline{Q_T})},||u||_{C(\overline{Q_T})}\right]$.
Our next step is to study partial derivatives of $\Phi(u)$. Simple application of the chain rule yields 
\begin{align}   
\partial_t\Phi(u)=&\Phi'(u)\partial_t u=\int_0^u\left(\int_0^x\varphi\left(\sfrac{y}{2}\right)\,dy\right)^2\,dx\,\partial_t u,\label{d2}\\
\nabla\Phi(u)=&\Phi'(u)\nabla u=\int_0^u\left(\int_0^x\varphi\left(\sfrac{y}{2}\right)\,dy\right)^2\,dx\,\nabla u,\label{d2_}\\
\partial_{x_ix_ j}\Phi(u)=&\Phi'(u)\partial_{x_ix_j}u+\Phi''(u)\partial_{x_i}u\partial_{x_j}u\nonumber\\
=&\int_0^u\left(\int_0^x\varphi\left(\sfrac{y}{2}\right)\,dy\right)^2\,dx\,\partial_{x_ix_j}u+\left(\int_0^u\varphi\left(\sfrac{x}{2}\right)\,dx\right)^2\,\partial_{x_i}u\partial_{x_j}u\text{ for } i,j\in 1:d.\label{d3}
\end{align}
In order to gain estimates for \eqref{d2}-\eqref{d3} in $C^{\left(\sfrac{\alpha}{2},\alpha\right)}(\overline{Q_T})$, we recall two properties of H\"older norms. Let $C$ be a closed subset of an open set $X\subset\R^d$ and let $\beta\in(0,1]$.  For all $g_1,h_1,h_2\in C^{(\beta)}(X)$, $g_2\in W^{1,\infty}(g_1(X))$, it holds that 
\begin{align}
 &||g_2\circ g_1||_{C^{(\beta)}(C)}\leq ||g_2||_{L^{\infty}(g_1(C))}+||\nabla g_2||_{L^{\infty}(g_1(C))}||g_1||_{C^{(\beta)}(C)},\label{rule1}\\
 &||h_2h_1||_{C^{(\beta)}(C)}\leq ||h_1||_{C^{(\beta)}(C)}||h_2||_{C^{(\beta)}(C)}.\label{rule2}
\end{align}
Using \eqref{rule1}-\eqref{rule2} and the properties of $\varphi$, we conclude from  \eqref{d2}-\eqref{d3} that
\begin{align}
 ||\Phi(u)||_{C^{\left(1+\sfrac{\alpha}{2},2+\alpha\right)}(\overline{Q_T})}\leq \Cl{CN}\label{est4}
\end{align}
for some $\Cr{CN}>0$ which depends only upon the parameters of the problem.
Indeed, for all $i,j\in 1:d$ and $k\in\left(0,\sfrac{1}{2}||u||_{C(\overline{Q_T})}\right]$, it holds, for instance, that
\begin{align}
&\left\|\Phi''(u)\partial_{x_i}u\partial_{x_j}u\right\|_{C^{\left(\sfrac{\alpha}{2},\alpha\right)}(\{k\leq u\leq 2k\})}\nonumber\\
=&\left\| \left(\int_0^u\varphi\left(\sfrac{x}{2}\right)\,dx\right)^2\partial_{x_i}u\partial_{x_j}u\right\|_{C^{\left(\sfrac{\alpha}{2},\alpha\right)}(\{k\leq u\leq 2k\})}\nonumber\\
\leq &\left\|\int_0^u\varphi\left(\sfrac{x}{2}\right)\,dx\right\|_{C^{\left(\sfrac{\alpha}{2},\alpha\right)}(\{k\leq u\leq 2k\})}^2\left\|\partial_{x_i}u\right\|_{C^{\left(\sfrac{\alpha}{2},\alpha\right)}(\{k\leq u\leq 2k\})}\left\|\partial_{x_i}u\right\|_{C^{\left(\sfrac{\alpha}{2},\alpha\right)}(\{k\leq u\leq 2k\})}\nonumber\\
\leq &\left\|\int_0^u\varphi\left(\sfrac{x}{2}\right)\,dx\right\|_{C^{\left(\sfrac{\alpha}{2},\alpha\right)}(\{k\leq u\leq 2k\})}^2\varphi^{-2}(k)\nonumber\\
\leq &\left(||u||_{C^{\left(\sfrac{\alpha}{2},\alpha\right)}(\overline{Q_T})}+2k\right)^2\nonumber\\
\leq &5||u||_{C^{\left(\sfrac{\alpha}{2},\alpha\right)}(\overline{Q_T})}^2.\nonumber
\end{align}
Similarly, for all $k\in\left(0,\sfrac{1}{2}||u||_{C(\overline{Q_T})}\right]$ it holds that
\begin{align}
 \left\|\Phi''(u)\partial_{x_i}u\partial_{x_j}u\right\|_{C^{\left(\sfrac{\alpha}{2},\alpha\right)}(\{-2k\leq u\leq -k\})}\leq 5||u||_{C^{\left(\sfrac{\alpha}{2},\alpha\right)}(\overline{Q_T})}^2.\nonumber
\end{align}
Consequently, we obtain that
\begin{align}
\left\|\Phi''(u)\partial_{x_i}u\partial_{x_j}u\right\|_{C^{\left(\sfrac{\alpha}{2},\alpha\right)}\left(\overline{\{u\neq 0\}}\right)}\leq5||u||_{C^{\left(\sfrac{\alpha}{2},\alpha\right)}(\overline{Q_T})}^2.\label{est6_}
\end{align}
Since $\Phi(u)\equiv0$ on $\{u=0\}$, \eqref{est6_} yields 
\begin{align}
\left\|\Phi''(u)\partial_{x_i}u\partial_{x_j}u\right\|_{C^{\left(\sfrac{\alpha}{2},\alpha\right)}(\overline{Q_T})}\leq 5||u||_{C^{\left(\sfrac{\alpha}{2},\alpha\right)}(\overline{Q_T})}^2.\label{est6}
\end{align}
Treating the remaining three terms on the right-hand sides of \eqref{d2}-\eqref{d3} in the same way, we obtain the estimate \eqref{est4}.

Now we are ready to produce the variable transformation $V$ with the desired properties. We define
\begin{align}
  &V:\left[-||u||_{C(\overline{Q_T})},||u||_{C(\overline{Q_T})}\right]\rightarrow\R,\ V(k):=\int_0^k(a^{-1}\Phi')(x)\,dx.\nonumber
\end{align}
Clearly, $V(0)=0$. Further, due to the properties of $\varphi$, particularly \eqref{sm}, it holds that
\begin{align}
 |V'(k)|=|(a^{-1}\Phi')(k)|=&a^{-1}(k)\left|\int_0^k\left(\int_0^x\varphi\left(\sfrac{y}{2}\right)\,dy\right)^2\,dx\right|\nonumber\\
 \leq&\sfrac{1}{3}|k|^3 
 a^{-1}(k)\varphi^2\left(\sfrac{k}{2}\right)\leq \Cl{C2}\text{ for all }k\in\left[-||u||_{C(\overline{Q_T})},||u||_{C(\overline{Q_T})}\right]\nonumber
\end{align}
for some constant $\Cr{C2}$ which depends only upon the parameters of the problem. Hence, the function $V$ is well-defined and belongs to $W^{1,\infty}\left[-||u||_{C(\overline{Q_T})},||u||_{C(\overline{Q_T})}\right]$. It is clear also that $V$ is strictly increasing. 
Let us now check that $\partial_t V(u)\in C^{\left(\sfrac{\alpha}{2},\alpha\right)}(\overline{Q_T})$. Again, it is sufficient to consider this function on the sets $\{k\leq u\leq 2k\}$ and $\{-2k\leq u\leq -k\}$ for arbitrary $k\in\left(0,\sfrac{1}{2}||u||_{C(\overline{Q_T})}\right]$ and prove that the H\"older norms are bounded by a constant which is independent of $k$. So let $k\in\left(0,\sfrac{1}{2}||u||_{C(\overline{Q_T})}\right]$. It holds due to the properties of $\Phi$, \eqref{rule1}-\eqref{rule2} and condition \eqref{sm} that
\begin{align}
&||\partial_t V(u)||_{C^{\left(\sfrac{\alpha}{2},\alpha\right)}(\{k\leq u\leq 2k\})}\nonumber\\
=&||(a^{-1}\Phi')(u)\partial_t u||_{C^{\left(\sfrac{\alpha}{2},\alpha\right)}(\{k\leq u\leq 2k\})}\nonumber\\
 \leq &||a^{-1}(u)||_{C^{\left(\sfrac{\alpha}{2},\alpha\right)}(\{k\leq u\leq 2k\})}||\Phi'(u)||_{C^{\left(\sfrac{\alpha}{2},\alpha\right)}(\{k\leq u\leq 2k\})}||\partial_t u||_{C^{\left(\sfrac{\alpha}{2},\alpha\right)}(\{k\leq u\leq 2k\})}\nonumber\\
 \leq &\Cr{CN}\left(a^{-2}(k)||a(u)||_{C^{\left(\sfrac{\alpha}{2},\alpha\right)}(\{k\leq u\leq 2k\})}+a^{-1}(k)\right)\varphi^{-1}(k)\nonumber\\
 \leq&\Cl{C4}a^{-2}(k)\varphi(k)\leq\Cl{C5}\label{est8}
\end{align}
for some constant $\Cr{C5}>0$, which once again depends only upon the parameters of the problem. 
Similar estimates hold on the sets $\{-2k\leq u\leq -k\}$.

Let us know go back to \eqref{PM1} and multiply it by $V'(u)$. We obtain after standard calculation that 
\begin{align}
 \partial_t V(u)=&\Phi'(u)\Delta u+(a^{-1}\Phi'f)(u)\nonumber\\
 =&\Delta\Phi(u)-\Phi''(u)\left|\nabla u\right|^2+(a^{-1}\Phi'f)(u)\text{ in } Q_T.\label{eqgood}
\end{align}
We define
\begin{align}
 &\bar{f}_u:=-\Phi''(u)\left|\nabla u\right|^2+(a^{-1}\Phi'f)(u)\text{ in } Q_T.\nonumber
\end{align}
Equation \eqref{eqgood} then reads:
\begin{align}
 \partial_t V(u)=\Delta\Phi(u)+\bar{f}_u\text{ in } Q_T.\label{eq6}
\end{align}
We already know from  estimates \eqref{est4} and \eqref{est8} that functions $\Delta\Phi$ and $\partial_t V(u)$ are in $C^{\left(\sfrac{\alpha}{2},\alpha\right)}(\overline{Q_T})$.  Therefore, equation \eqref{est6} holds in $C^{\left(\sfrac{\alpha}{2},\alpha\right)}(\overline{Q_T})$.
Finally, since $\Phi'(u)=(aV')(u)$, we can rewrite \eqref{eq6} as
\begin{align}
 \partial_t V(u)=\nabla\cdot(a(u)\nabla V(u))+\bar{f}_u\text{ in } Q_T.\label{eq4}
\end{align}
For the new variable $v:=V(u)$, \eqref{eq4} takes the form \eqref{eq3}:
\begin{align}
 \partial_t v=\nabla\cdot(a\circ U(v)\nabla v)+\bar{f}_u\text{ in } Q_T,\nonumber
\end{align}
where $U$ is the inverse of the function $V$. This finishes the proof of {\it Theorem~\ref{maint}}.\\\\
\qed
\section{Example}\label{Example}
Let us consider for $m>1$ the homogeneous porous medium equation
\begin{align}
 &\partial_t u=m u^{1-\sfrac{1}{m}}\Delta u\text{ in }[0,T]\times \overline{B_R}=:\overline{Q_T},\label{PMB}
\end{align}
where $T>0$ and $B_R$ is the open $d$-ball of an arbitrary radius $R>0$ centred at the origin. One of the solutions of \eqref{PMB} is the function
\begin{align}
 u(t,x)=B(t+1,x)\text{ for all }t\in[0,T],\ x\in B_R,\nonumber
\end{align}
where $B$ is the well-known Barenblatt solution \cite{Barenblatt}
\begin{align}
 B(t,x)=t^{-\frac{md}{d(m-1)+2}}\left(C-\frac{m-1}{2m(d(m-1)+2)}|x|^2t^{-\frac{2}{d(m-1)+2}}\right)_+^{\frac{m}{m-1}}\nonumber
\end{align}
for arbitrary constant $C>0$.
Let $\alpha\in(0,1)$ be arbitrary. We set
\begin{align}
 \Phi(k):=\frac{m(1+\alpha)}{2+\alpha}k^{\left(1-\sfrac{1}{m}\right)(2+\alpha)}\text{ for all }k\in\R^+_0.\nonumber
\end{align}
Clearly, $\Phi(u)\in C^{(2+\alpha)}(\overline{Q_T})$. Following the construction from the proof of {\it Theorem~\ref{maint}}, we introduce the variable transformation 
\begin{align}
 v:=V(u):=\int_0^u(\sfrac{1}{m}u^{\sfrac{1}{m}-1}\Phi')(x)\,dx=u^{{\left(1-\sfrac{1}{m}\right)(1+\alpha)}}\nonumber
\end{align}
and the new reaction term 
\begin{align}
 &\bar{f}_u:=-\Phi''(u)\left|\nabla u\right|^2=-(1+\alpha)(m-1)({\left(1-\sfrac{1}{m}\right)(2+\alpha)}-1)u^{{\left(1-\sfrac{1}{m}\right)(2+\alpha)}-2}|\nabla u|^2.\nonumber
\end{align}
Again, it is clear that $V(u)\in C^{(1+\alpha)}(\overline{Q_T})$. Thus, we obtain that $v$ is a classical solution of a porous medium equation: 
\begin{align}
 \partial_t v=m\nabla\cdot\left(v^{\frac{1}{1+\alpha}}\nabla v\right)+\bar{f}_u\text{ in }C^{(\alpha)}(\overline{Q_T}).\nonumber
\end{align}
\section*{Acknowledgement}
The author gratefully acknowledges the financial support from the Wolfgang Pauli Institute in Vienna during her stay there through the summer and autumn of 2014.
\printbibliography
\end{document}